\documentclass[preprint,authoryear,5p,times,twocolumn]{elsarticle}

\usepackage{graphics}
\usepackage{amssymb}
\usepackage{multirow}

\begin{document}
\bibliographystyle{model2-names}

\begin{frontmatter}
\title{Strategies in tower solar power plant optimization.}

\author[label1]{A. Ramos}
\author[label2]{F. Ramos}
\address[label1]{Centre de Physique Th\'eorique
    \fnref{label3}, CNRS Luminy, Case 907, F-13288 Marseille Cedex 9,
    France. {\tt <alberto.ramos@cpt.univ-mrs.fr>}}
\address[label2]{Nevada Software Informatica S. L. {\tt <francisco.ramos@nspoc.com>}}
\fntext[label3]{CPT is research unit
UMR 6207 of the CNRS and of Aix-Marseille univ. and Univ. Sud
Toulon-Var.}

\begin{abstract}
  A method for optimizing a central receiver solar thermal electric
  power plant is studied. We parametrize the plant design as a
  function of eleven design variables and reduce the problem of finding
  optimal designs to the numerical problem of finding the minimum of a
  function of several variables. This minimization problem is attacked
  with different algorithms both local and global in nature. We find
  that all algorithms find the same minimum of the objective
  function. The 
  performance of each of the algorithms and the resulting designs are
  studied for two typical cases.

  We describe a method to evaluate the impact of design variables
  in the plant performance. This method will
  tell us what variables are key to the optimal plant design and which
  ones are less important. This information can be used to further
  improve the plant design and to accelerate the optimization
  procedure. 
\end{abstract}


\end{frontmatter}


\section{Introduction}
\label{sc:intro}

One of the main tasks in the conversion of solar energy into
electricity by solar power plants is to work out an optimized plant
design. In this type of plant, the energy collector subsystems
(heliostats, field receivers) represent a very important part of the
cost break-down structure. Therefore, the use of
detailed computer programs is of great interest in order to optimize
the plant design.

The two main conceptual ingredients for a solar plant
optimization code are:
\begin{enumerate}
\item Reduction of the plant design to the value of certain
  \emph{design variables}. 
\item An optimization criteria: This means having a function that
  computes the objective quantity (i.e. total 
  annual power output, cost per produced power, etc\dots) as a
  function of the design variables. In general we will use the cost of
  the energy produced by the plant as the optimization criteria, and
  therefore use the terms \emph{optimize} and \emph{minimize} as
  interchangeable. 
\end{enumerate}

After the plant design is completed it becomes crucial to understand
what role each 
of the design variables play in the optimal design. Some variables can
be slightly moved away from the optimal value without impacting the
plant performance while others can not be changed without a
severe impact in the plant performance. We will give a precise
definition of a quantity (we call it \emph{uncertainty}) that will
measure the importance of each variable in a plant design. We
will give a precise mathematical definition of this \emph{uncertainty}
associated with a variable, and show how to compute these quantities
for a general plant design. 

The paper is organized as follows. In section~\ref{sec:desc} we will
explain our choice of plant design 
variables. Sections~\ref{sec:input} and~\ref{sec:per} describe
additional information needed to perform the
optimization. We perform the numerical optimization using different
algorithms, both global and
local. Section~\ref{sec:opt} describes in detail the optimization
procedure and our choice of three different algorithms: first we use a
fast local optimizer especially designed to solve this problem, second
we use complex optimization library that include a Monte-Carlo search,
with the potential ability to jump over function barriers and find
global optima. Finally we use a genetic algorithm, generally used for
difficult optimizations and problems with a complex fitness
landscape. 

Section~\ref{sec:prop} attack the very important problem of
determining what role each of the design variables have in the plant
performance. As we
will see, the Hessian of the objective function evaluated at the
minimum will provide information on how flat each of the directions in
the minimum are. This information will not only give us the
uncertainty associated with each design variable, but also help us
deciding if
the function have several local optima and if the three algorithms
have found the same design.

Finally, in section~\ref{sec:res} we apply these optimization
algorithms and analysis tools to two typical plant designs, and in
section~\ref{sec:conc} we summarize our main results and present some
perspectives for further studies. 

\section{Description of plant design variables}
\label{sec:desc}

In choosing the variables that determine the layout of our
plant, one has to take into account that an optimization procedure
requires multiple evaluations of the objective function. Thus a
proper choice of design variables must have the CPU cost
of the optimization in mind.

In this section we will present our choice for plant design
variables. This choice has been made with the following things in mind
\begin{itemize}
\item A plant should be circular-like to minimise the blocking and shadowing
  effect. 
\item Consecutive files of heliostats should be allocated in
  \emph{radial staggered} position to minimise blocking and shadowing
  effect, and keep a compact field.
\end{itemize}

\begin{figure}
  \centering
  \includegraphics[width=8cm]{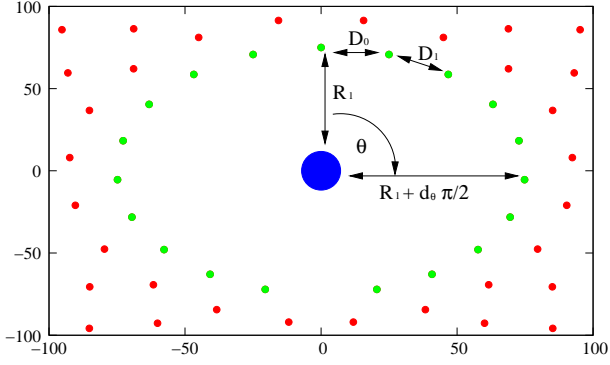}
  \caption{An example of the first two lines of heliostats in a
    field. The azimuthal angle starts in the north with the value
    $\theta=0$. In green we have the heliostats that belong to the
    first line. The distance of this line to the tower is given by
    $R_1$. The distance of all the lines $R_n$ is determined trough
    the recursion relation of Eq.~\ref{eq:recd}. We see that this
    distance dependes on the azimuthal coordinate through the variable
    $d_\theta$ (Eq.~\ref{eq:rcorr}). The azimuthal distance between
    heliostats of the first 
    line ($D_0^1, D_1^1,\dots$) are determined through the recursive
    relation Eq.~\ref{eq:reca}. The azimuthal position of the second
    and subsequent lines of heliostats is determined by positioning
    them at the average azimuthal position of the heliostats of the
    previous line (radial staggered) until a transition line is
    placed. 
  }
  \label{fig:collector}
\end{figure}

\subsection{Collector field variables}

The distance of the $n^{th}$ heliostat line to the tower ($R_n$) is
given by the linear recursion formulae
\begin{equation}
  \label{eq:recd}
  R_n = a_0 + (a_1 + 1)R_{n-1}
\end{equation}
that defines the two variables
\begin{enumerate}
\item[V01] $a_0: $ Initial spacing between heliostat rows $[m]$.
\item[V02] $a_1: $ Increasing spacing between heliostat rows.
\end{enumerate}

The recursion relation of Eq.~\ref{eq:recd} is solved by the condition
$R_0 = R_{\sc base}$ that specifies the distance between the tower
and the first heliostat line. This space is usually used for
operational purposes, like the administrative building, roads,
etc\dots There is a minimal distance between heliostats 
$R_{\sc min}$ related with the size of the heliostats and some
practical needs. $R_n$ is set to the maximum between the
value given by the recursion relation of Eq.~\ref{eq:recd} and $R_{\sc
min}$. 

The shape of the lines of heliostats of an optimal field does not have
to be perfectly circular. We consider the possibility of non circular
shapes by adding to the previous radial distance an increment that
dependes on the azimuthal position of the heliostat (we will call it
$\theta$). 
\begin{equation}
  \label{eq:rcorr}
  \Delta R = \left\{
  \begin{array}{lcl}
    d_\theta \theta & \textrm{for} & 0 < \theta \le \pi \\
    d_\theta (2\pi - \theta) & \textrm{for} & \pi < \theta < 2\pi \\
  \end{array}
\right. .
\end{equation}
This defines another collector variable
\begin{enumerate}
\item[V03] $d_\theta: $ Correction of the radial distance with the
  azimuthal position.
\end{enumerate}
Intuitively when $d_\theta > 0$ the lines of heliostats will be closer
to the tower in the north part of the field, whereas for $d_\theta <
0$ heliostats will be closer to the tower in the south. Usually
optimal designs are more compact in the south (and hence $d_\theta <
0$) where there is less blocking and shadowing effect between heliostats.  

A proper optimal layout should not only find the optimum radial
spacing between heliostats of the same line, but also determine the
azimuthal distance between them. We will use a 
similar technique than the one used for the radial distance. We
need to determine how this azimuthal distance depends on both the
azimuthal angle and the radial distance. If we call $D_0^1$ the
azimuthal distance for the heliostat that is just in 
the north (zero azimuth) in the first line of heliostats, and 
number the heliostats in the same line with the index $\alpha=0,1,\dots$, 
the azimuthal distance as a function of the azimuthal angle $\theta$
is given by
\begin{equation}
  \label{eq:reca}
  D_\alpha^1 = D_{\alpha-1}^1 + e_\theta\theta
\end{equation}
that, defines an additional variable
\begin{enumerate}
\item[V04] $e_\theta: $ Variation of the azimuthal distance with the
  azimuthal angle. 
\end{enumerate}
To start this recursion relation one should give $D_0^1$ as input. We
will comment about this later.

Variable ([V04]) determine the azimuthal position only of the
first line of heliostats. The remaining lines of heliostats are
situated at the radial distance determined by the variables
([V01--V03]) and with radial staggered positions. This means that their
azimuthal angle is the average of the azimuthal position of the two
heliostats in front of it (see Fig.~\ref{fig:collector}). 

\begin{figure}[h!]
  \centering
  \includegraphics[width=7cm]{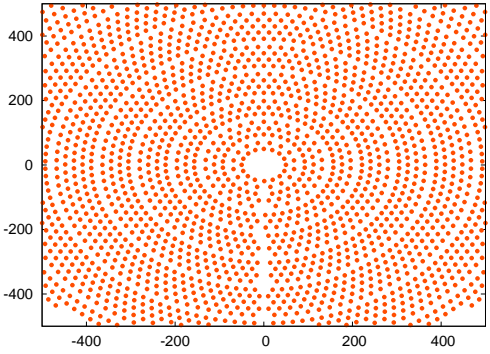}
  \caption{Close look to the first transition lines of the field
    layout of an optimal 3000 
    heliostat plant. We can see the transition lines, the number of
    heliostats in each these transition lines are 2-3-3-4-4-5-... (see
    text for more details). This optimal layout correspond to the output of 
    the NSPOC code (see Sections~\ref{sec:opt} and~\ref{sec:res} for
    more details). }
\label{fig:c3000}
\end{figure}

This rule of formation will increase azimuthal spacing between
heliostats very fast. So it is convenient, after a certain amount of
lines, to give an extra space and restart again the rule of formation
forgetting the previous heliostats. This extra space increases with
the distance to the tower. These re-starting lines are called
\emph{transition lines}, and heliostats between two transition lines
are said to belong to the same \emph{group}. The first line of
heliostats of each group is at a distance from the
previous heliostat line given by 
\begin{equation}
  \label{eq:trans}
  \Delta = (1 + a_0 + a_1R)\delta + \epsilon
\end{equation}
where $R$ is the last distance between lines given by
Eq.~\ref{eq:recd}. Eq.~\ref{eq:trans} defines two more optimization
variables 
\begin{enumerate}
\item[V05] $\epsilon: $ Extra distance for transition lines.
\item[V06] $\delta: $ Increment of distance in transition lines with
  distance from the tower.
\end{enumerate}
The variable $\epsilon$ determines the extra spacing in a
transition line to avoid large blocking and shadowing effects. Variable
$\delta$ sets an extra space that is proportional to the radial
distance of the transition line. In Fig.~\ref{fig:c3000} we can see
the transition lines of a 3000 heliostats plant.  

After a transition line we have to determine the azimuthal distance
again of the next group of heliostats. This is done with the same 
recursion relation Eq.~\ref{eq:reca}. 
\begin{equation}
  D_\alpha^a = D_{\alpha-1}^a + e_\theta\theta
\end{equation}
where the index $a=1,2...$ label groups of heliostats. As was
commented before these recursion relations need $D_0^a$ as an initial
condition. These quantities are determined through
\begin{equation}
  D_0^a = (b + 1)D_0^{a-1}
\end{equation}
that essentially allow the azimuthal distance between heliostats in
the north to increase/decrease with the radial distance of the
group. This last recursion relation needs $D_0^1$ as an input. This is
another optimization parameter completing 
the 8 parameters needed to optimize a field layout
\begin{enumerate}
\item[V07] $b: $ Azimuthal distance dependence with the radial
  distance of the group.
\item[V08] $D_0^1: $ Initial azimuthal distance in the first line of
  heliostats of the plant.
\end{enumerate}
It is understood that there is a minimal azimuthal distance $D_{\sc
  min}$ between 
heliostats that is given by the size of the heliostat plus some
arbitrary distance needed for operational purposes.

The only additional information needed is the number of heliostat
lines of
each group. These are determined by trial and error. Reasonable results
are obtained with an increasing sequence like
$2-3-3-4-4-5-5-6-6-7-7...$ (see Fig.~\ref{fig:c3000} for an
example). The conclusions of this paper are 
unchanged by the details like how many heliostats each group has, one
only needs to fine tune this numbers when looking for the final design
of a plant, and this can be easily done in an automated way. 

It is important to remark here that the set of variables ([V01--V08])
are compatible with much more conventional field designs. For example
having no transition lines is easily achieved by adjusting $\epsilon$
and $\delta$. A
constant azimuthal distance between heliostats by setting
$b=e_\theta=0$, and a constant separation between lines corresponds to the
choice $a_1=0$. Finally $d_\theta = 0$ give perfectly circular heliostat
lines. 

The existence of some of the variables allow more complex designs that
we think can improve the plant performance, but this complex design
is not imposed. The optimization process will choose between the
different options.

Following these rules one generates a general field layout. We chose
to fix the number of heliostats of our plant design ($N_{\rm
  hel}$). In order to choose these heliostats we pick up the ones
which contribute with more power to the receiver. 
\begin{figure}[h!]
  \centering
  \includegraphics[width=7cm]{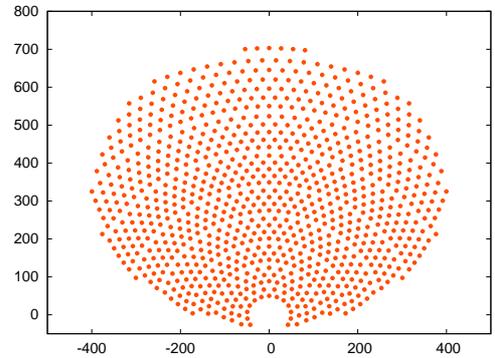}
  \caption{Field layout of an optimal 900
    heliostat plant. This optimal layout correspond to the output of
    the NSPOC code (see Sections~\ref{sec:opt} and~\ref{sec:res} for
    more details). During the optimization more heliostats are
    generated but not selected for the final layout due to their worse
    performance compared to the selected ones (see
    Fig.~\ref{fig:n900-no}).} 
\label{fig:n900}
\end{figure}

\begin{figure}[h!]
  \centering
  \includegraphics[width=7cm]{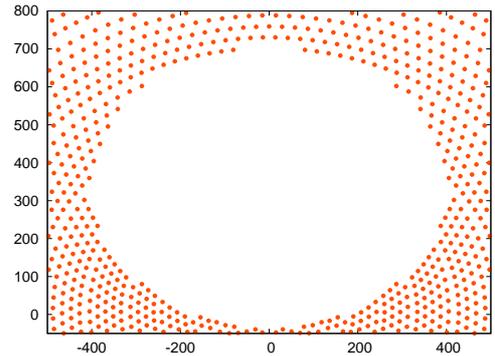}
  \caption{Heliostats not selected for the final layout due to their
    worse performance. This image shows a subset of all the generated
    heliostats during the optimization (see text and
    Fig.~\ref{fig:n900} for more details.)}  
\label{fig:n900-no}
\end{figure}

In Fig.~\ref{fig:n900} and Fig.~\ref{fig:n900-no} we can see how this
procedure works in practice 
for a 900 heliostat plant. During the optimization process, we
generate a field layout with far more than 900 heliostats. To draw the
final layout we simply pick the ``best'' heliostats of the field (those
shown in Fig.~\ref{fig:n900}), but more heliostats are generated and
rejected due to their worse performance (see Fig.~\ref{fig:n900-no}).

\begin{figure}
  \centering
  \includegraphics[width=7cm]{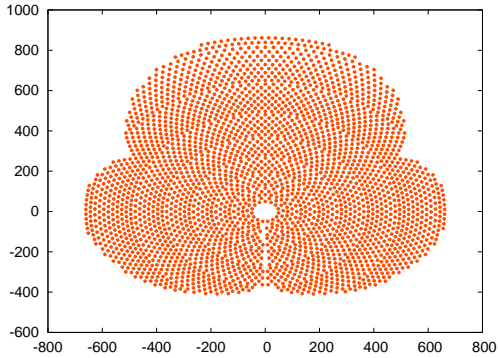}
  \caption{Field layout of a plant with a tower with three
    receivers. One of them points to the north, and the others at
    120$^{o}$. As can be seen the presence of two additional receivers
    increases dramatically the interception efficiency of heliostats at
    the east and west, making them to contribute with more power to
    the receiver than other heliostats in the north but more
    distant. This effect make the optimal field layout to have a
    clover shape.} 
  \label{fig:3cavity}
\end{figure}

Note that this rule to select the heliostats does not necessarily
impose a circular-like field layout. Simply this layout turns out to
be optimal for the problems we are studying here. As an example of a
non-circular optimal field, 
we show in Fig.~\ref{fig:3cavity} the field layout of a more
complicated plant design with three cavity receivers in the
tower (see~\cite{Ramos:2009aa} for more details). As can be seen in
Fig.~\ref{fig:3cavity} having three cavity receivers dramatically improves
the interception efficiency of some heliostats in the east and west
part of the field, making this heliostats to perform better that
heliostats that are far in the north. In this case
the optimal field layout does not have a circular shape, but clover-like.

\subsection{Receiver variables}

We can design either a north field plant or a circular plant. In the
first case we will consider that the receiver consists on an circular
aperture\footnote{We have considered more complex aperture shapes,
  like elliptical or rectangular. Although some of this designs
  perform better than a circular aperture, the conclusions of
  the paper remains completely unaltered, and keep the discussion
  simpler.} pointing to the north that lets the reflected solar rays
to enter. This receiver is characterized by the following variables 
\begin{enumerate}
\item[V09] $h_T: $ Tower height $[m]$.
\item[V10] $r:   $ Aperture radius $[rad]$.
\item[V11] $e_L: $ Aperture inclination $[rad]$.
\end{enumerate}

In the second case the receiver consists in a cylinder that can absorb
radiation coming from all directions. In this case the receiver is
characterized by its position and size\footnote{In real world designs
  one should worry about the maximum power density absorbed by the
  receiver, since this is constrained by material properties. This can
  (and should) be included in the optimization process, but we will not
  address this problem here.}, parametrized by the
following variables 
\begin{enumerate}
\item[V09] $h_t: $ Tower height $[m]$.
\item[V10] $r:   $ Receiver radius $[rad]$.
\item[V11] $h_r: $ Receiver height $[m]$.
\end{enumerate}

\section{Additional input: minimization parameters}
\label{sec:input}

To compute the full output of a solar power plant, we need
some additional input related with the heliostat characteristics and
plant location. These parameters are not treated as variables in the
minimization process. 

These extra input are seven parameters that define the heliostats
characteristics and plant location, plus insolation data from the plant
location. 

The study on how the plant performance depends on these parameters
(what we call \emph{parametric analysis}) is definitively a
very interesting subject that can answer questions like What is the
proper heliostat size? How much influence the optical quality of an
heliostat the plant performance?, etc... Nevertheless these questions
are beyond the scope of the present work and need further
investigations~\citep{Ramos:2011aa}.

\subsection{Heliostat characteristics}
\label{sec:hel}

An heliostat is characterised by its geometry and its optical
properties. All heliostats are assumed to be rectangular,
with the focal equal to the slant range and made of spherical
facets. 

We have a total of four parameters to describe these
properties of heliostats. All these parameters can be different for
different groups of heliostats (what we call \emph{mixed fields}):

\begin{enumerate}
\item [P01] $\sigma_h: $ Heliostat optical error.
\item [P02] $L_h: $ Horizontal length of the heliostat.
\item [P03] $L_v: $ Vertical length of the heliostat.
\end{enumerate}

\subsection{Geographic characteristics}

The location and local ground characteristics of the solar power plant
are taken into account in the following parameters
\begin{enumerate}
\item [P05] $\phi: $ Latitude of the plant position.
\item [P06] $m_N: $ Terrain north-south slope.
\item [P07] $m_W: $ Terrain east-west slope.
\end{enumerate}

\subsection{Insolation and ambiente temperatures}

Direct insolation in clear days and ground level temperatures at day
hours the 21$^{\textrm{st}}$ of each month are part of the input data.

Insolation ratios for each month are also input data. This ratio is
defined as the relation between the solar energy received in a  month
and the solar energy received if all days were clear.

\section{Performance models}
\label{sec:per}

We will not detail all the models used to compute the
performance of the power plant, only to say some words and give the
appropriate references. 

The interception efficiency is computed with an optical model
described in~\citep{Kierz:1980aa,Kierz:1980xs}. 

Atmospheric losses are estimated by the model described
in~\citep{Biggs:1976}.  

Receiver and cycle efficiency is estimated by taking into account the total
input power through the aperture, the aperture size, the aperture
inclination and the ambient temperature. This estimation is done by
using a polynomial function whose coefficients have been determined by
fitting data. 

Finally the total daily and yearly energies are computed by
integrating the hourly power output. We compute total daily and yearly
energies for clean day, and mean cloudy basis. 

\section{Optimization procedure}
\label{sec:opt}

Within the framework described in the preceding sections we are in the
position to start an optimization procedure. By using the performance
models described above we can estimate the total yearly energy output,
and with this our optimization citeria, the price of the produced
energy, as a function of our eleven design variables ([V01--V11]), seven
design parameters ([P01--P07]) and insolation data. 
\begin{equation}
  E(v_i; p_\alpha, I)
\end{equation}
where $v_i$ ($i=1,\dots,11$) represents the optimization variables,
$p_\alpha$ ($\alpha=1,\dots,7$) are the design parameters, and $I$ is
the insolation data. 

One only needs to find the values of $v_i$ that makes $E(v_i;
p_\alpha, I)$ maximum. Being $E(v_i; p_\alpha, I)$ a non linear
function one should be concerned about the existence of local optima. 

To have a full control over the optimization process we have made some
tests on typical plant configuration by using several state of the art
algorithms. First we will use our own local optimization
algorithm~\citep{Paco:nspoc}, 
designed with the particularities of solar plant optimization in
mind. This is a fast algorithm that find the closest local
optimum. Second we will use a mixture of Monte Carlo and
conjugate gradient like algorithms coded in the CERN \texttt{MINUIT}
library~\citep{James:1975dr}. Finally we have specifically coded a
genetic algorithm~\citep{Ramos:ga} commonly used to solve difficult
optimization problems with many local optima. 

Algorithms capable of finding global optima are usually much slower
than local algorithms. The porpoise of using global searchers is to
find out if our preferred local algorithm finds the same optima as
global algorithms.

Now we will describe in detail our optimization algorithms.

\subsection{\texttt{NSPOC} algorithm}

Our preferred local algorithm is a variant of Powell's
algorithm~\citep{Powell:1964aa}. 

Starting from an approximate value of
the optimal variables $v_i^{(0)}$ our algorithm performs a line search
along the first direction ($v_1$) until a
optimum is found with some initial precision $h$. Then the algorithm
proceeds with the second direction, and so on, until the 11 variables
have been explored. This process needs to be repeated until all the
variables are fixed at their optimum for a full cycle. Finally one can
increase the precision of the optima position by repeating the 
process with an increased precision (i.e. setting $h\rightarrow h/2$).
\begin{figure}
  \centering
  \includegraphics[width=8cm]{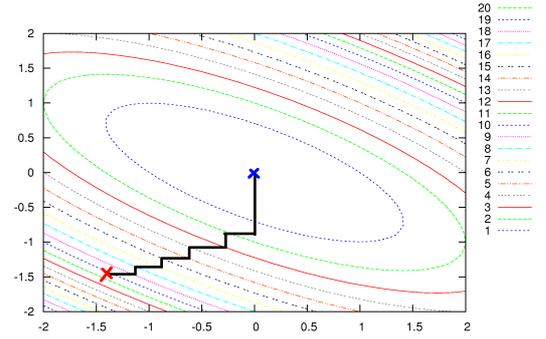}
  \caption{"Zig-zagging" of the \texttt{NSPOC} algorithm to find the
    minimum in an example two dimensional function}
\label{fig:cnt}
\end{figure}

Following this algorithm we will ``zig-zag'' until the optima is
found (see Fig.~\ref{fig:cnt}). One may naively think that this
zig-zag results in a very slow convergence, but this turns out not to
be the case. Exploring the dimensions one by one instead of
approaching the minima along a ``direct'' route has the advantage that
one does not have to fully evaluate the objective function at each
step. For example by changing the aperture size or inclination one
does not need to recompute the blocking and shadowing effects, or the
atmospheric efficiency. 

Moreover algorithms that approach the minima along a ``direct'' route
needs additional information from the objective function in order to
find the proper direction in which the optima is located. This extra
information usually comes from the derivative of the
objective function (i.e. the conjugate gradient
method~\citep{Hestenes:1952aa}). We do not have the possibility of
computing analytically the derivative of the objective function, and a
numerical evaluation of the derivative severely worsen the convergence
speed of these algorithms. There also exists methods that tries to avoid the
`''zig-zagging'' without using the derivative (like the original
Powell method), but in practice this methods can easily end up
searching the optima in a lower dimensional subspace. 

In summary, cause the re-evaluation of the objective function is much
simpler when only one variable has changed, and because more
complicated ingredients need additional knowledge of the objective
function to work (like the gradient) our simple variation of Powell's
method turns out to be a fast and robust method for solar power plant
optimization.

\subsection{\texttt{MINUIT} algorithm}

The \texttt{MINUIT} library has been widely used in high energy
physics as well as in other fields with literally thousands of papers
based in its results. It is 
considered a robust minimizer. This is the reason we 
choose to use it here. We refer the reader to the \texttt{MINUIT}
reference~\citep{James:1975dr} for further details. Here we will only
comment that the \texttt{MINUIT} library has several algorithms
implemented. 

We first use the \texttt{SEEK} optimizer. This is a Monte Carlo search
that has potentially the ability to jump over function barriers to
find a better global optima. The algorithm includes a metropolis
step~\citep{Metropolis:1953aa} which moves to a new position
$v_i^{(new)}$ from an old one ($v_i^{(old)}$) with a
probability 
\begin{equation}
 \mathcal{P}\left(\nu_i^{(old)}\rightarrow \nu_i^{(new)}\right) 
 \propto \exp\left\{-\frac{E(v_i^{(new)}; p_\alpha,
    I)}{E(v_i^{(old)}; p_\alpha, I)}\right\} 
\end{equation}

We refine this optima search with the use of the \texttt{MINUIT}
optimizer \texttt{MIGRAD}, that is consider \texttt{MINUIT} best
local optimizer. It is a quasi-Newton method~\citep{Davidon:1991aa} with
inexact line search  and a stable metric updating scheme. We decided
to use the strategy that make less use of the numerical estimates of
the gradient of the objective function at the price of being slower.

Regardless the potential abilities of this combination of
\texttt{MINUIT} algorithms it is fair to say that about half of the
\texttt{MINUIT} hackers believes that the ability of finding global
minima are small in practical situations.

\subsection{Genetic algorithm}

To be sure that our local optimizer and \texttt{MINUIT} are not
falling in local optima, we have also implemented a genetic
algorithm ~\citep{Ramos:ga}\footnote{The \texttt{FORTRAN} code of the genetic
  library is free software under the GPL license.}.

Genetic algorithms~\citep{Holland:1975aa} mimic the process of natural
evolution by 
creating a ``population''  and implementing natural-inspired
mechanisms as crossover, mutation and natural selection to optimize an
objective function.

The main ingredients of a genetic algorithm are
\begin{enumerate}
\item A population of organisms whose fitness is given by the
  objective function. 
\item A crossover process, by which two members of the population
  (``parents'') give rise to two different members (``childs'').
\item A mutation process, by which organisms randomly change.
\item A selection process by which the members of the population are
  chosen for later crossover and/or mutation based in its fitness. 
\end{enumerate}

\subsubsection{Population}

The members of our population have 11 ``genes'' ($g_i,\,
i=1,\dots,11$) whose values are  the 11 
variables that code the plant design. In this case these variables are
stored simply as a vector of 11 real numbers. 

Members of the population have a ``fitness'' that is given as the
value of the objective function evaluated for its genetic content.
\begin{equation}
  \textrm{Fitness} = E(g_i; p_\alpha, I)
\end{equation}

We have repeated the optimization procedure with different population
sizes in the range $N_{tot} = 30-100$. 

\subsubsection{Crossover and mutation}

The population is paired up, and each pair crosses with a probability
$p_c$ (typical values in our runs are $p_c\approx 0.05$). When a
crossover occurs each of the genes of the two members of the
population $g_i^{(1)}$ and $g_{i}^{(2)}$ produce the offspring with
genetic content $g_i^{'(1)}$ and $g_{i}^{'(2)}$ given by
\begin{eqnarray}
  g_i^{'(1)} &=& r_ig_i^{(1)} + (1-r_i)g_i^{(2)} \\
  g_i^{'(1)} &=& q_ig_i^{(1)} + (1-q_i)g_i^{(2)} 
\end{eqnarray}
where $r_i, q_i$ are samples of a normal distribution $\mathcal
N(0.5,0.5)$ with mean $1/2$ and standard deviation $1/2$. It is important
here to remark that this choice allow interpolation as well as
extrapolation between parent's genetic content, avoiding a fast
``false convergence'' (see for example~\cite{Allanach:2004my}).

After this process each gene of each member of the population is
mutated with a probability $p_m$ (typical values in our runs are
$p_m=0.1$). In this mutation process the gene value is multiplied by
$(1+r)$ where $r\sim \mathcal N(0,1)$ is a random sample of a standard
normal distribution.

\subsubsection{Selection}

If a member is the product of a crossover or has mutated its fitness
(value of the target function) needs to be re-computed. 

The ``best'' $N_{elite}$ of the total population $N_{tot}$ organisms
are retained from one generation to 
the next without change. For the rest we use the well known roulette-wheel
selection process (see for example~\citep{Goldberg:1989aa} and
references therein) in which the members of the
population are ordered 
by its fitness and accumulated normalized fitness values are computed
(the accumulated fitness value is the sum of the fitness of all the
individuals better fitted than itself). After we choose
$N_{tot}-N_{elite}$ individuals by drawing uniform random numbers and
taken the first individual whose accumulated normalized fitness value
is greater than these random samples.

The idea behind this process is that if an organism has a fit that is
twice than other organism's fitness, it is twice as probable for this
organism to ``survive''.

This members together with the $N_{elite}$ members form the ``next
generation'' that is subject to the same process again. Typically
$N_{elite}$ is chosen 5-10. 

\subsubsection{Stop criteria}

The stop criteria is always a subtle matter in a genetic algorithm. We
choose to iterate the process a fixed number of generations (in
practice around 200), and examine if the fittest organism has changed
in the last 20 generations. If this is the case we proceed to iterate
for another 200 generations.

In all situations we are sure that when we decide to stop the
optimization process, the fitness of the elite ($N_{elite}$)
organism has not significantly changed in at most 20 generations. 

We also routinely perform tests during the run to ensure that there
exists diversity among the organisms to avoid a false fast
convergence. The existence of diversity is 
controlled by computing the variance of each genes in all
generations. 

\section{Impact of the variables in the plant performance}
\label{sec:prop}

\begin{figure}[h]
  \centering
  \includegraphics[width=9cm]{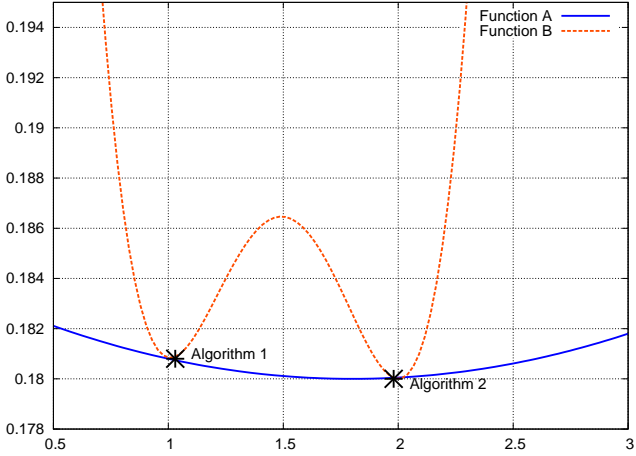}
  \caption{One dimensional example of the minimization of a complex
    function. Algorithm 1 outputs that the minimum is at $x=1.03$ with
    a value of the objective function equal to $1.1808$. Algorithm 2
    says that the minimum is located at $x=1.98$ with a value of the
    objective function equal to $1.18001$. These results are
    compatible with both function types A and B, but the physical
    interpretation is very different in both cases. A detailed study
    of the second derivative of the objective function will help us in
    discerning both cases (see text for more details).}
  \label{fig:multi}
\end{figure}

Once one finds the value of the 
eleven variables that makes the objective function optimal,
it is interesting to analyze the properties of this optimum. 
The basic idea is better understood in a one
dimensional case (see Fig.~\ref{fig:multi}).
When two algorithms outputs the position of the minimum and the value
of the objective function in the minimum they will never give exactly
the same answer. When the value of the objective function is very
similar, the difference in the position of the minimum can be due to
the fact that the minimum is very wide (function A of
Fig.~\ref{fig:multi}), or that the function has several local minimum
(function B of Fig.~\ref{fig:multi}). These cases can be easily
discerned by looking at the second derivative at the
minimum. Furthermore, the same analysis 
also provides information on how the value of different variables
affect the plant performance, helping us to detect ``irrelevant''
variables that we could have included in the design, and to focus in
variables that are key for the plant performance. This analysis for
the general case of an objective function of several variables will be
explained in this section and in the~\ref{ap:hess}. 

\begin{figure}
  \centering
  \includegraphics[width=7cm]{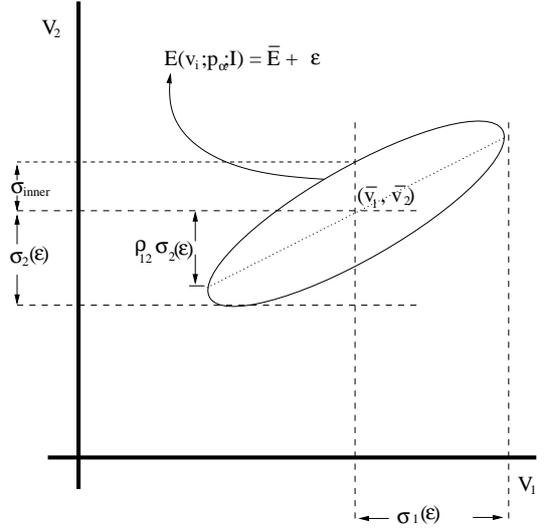}
  \caption{Contour plot of the function $E(\nu_i; p_\alpha, I)$ in the
    $\nu_1,\nu_2$ plane. The minimum of the function is located at
    $\nu_1 = \bar\nu_1$ and $\nu_2 = \bar\nu_2$ (at the center of the
    ellipse), and the value of the function in this point is $\bar
    E$. The ellipse corresponds to the points that having the values
    of $\nu_1$ and $\nu_2$ that correspond to the point in the plane,
    have a minimum value of the objective function equal to $\bar E +
    \varepsilon$. The size of the
    ellipse defines the quantities $\sigma_1(\varepsilon)$ and
    $\sigma_2(\varepsilon)$. The inclination define the
    ``correlation'' between the two variables $\rho_{12}$. This
    correlation can also be estimated from the difference between
    $\sigma_{\rm inner}$ and $\sigma(\varepsilon)$.}
  \label{fig:errorelipse}
\end{figure}

We will assume that the optimal value of the objective function is
achieved at $\bar \nu_i$ and that the value of the objective function
at this point is $E(\bar \nu_i; p_\alpha, I) = \bar E$. If we take the
value of one variable away 
from its optimal value by an amount $\sigma(\varepsilon)$ the value of
the objective function will increase. But some of the rise in the
objective function can be compensated by tuning the other
variables. We define the \emph{uncertainty of the variable $\nu_i$}
(and use the symbol $\sigma_i(\varepsilon)$) with the condition that 
the minimum of the objective function with the $i-th$ variable fixed
at the value $\nu_i = \bar\nu_i + \sigma_i(\varepsilon)$ is equal to
$\bar E + \varepsilon$ (see Fig.~\ref{fig:errorelipse}).

If $\varepsilon$ is chosen small, a large value of
$\sigma_a(\varepsilon)$ means that the variable $\nu_a$ plays no role
in the plant performance. It does not matter the value we choose for 
$\nu_a$ the objective function will not deviate from its optimal
value. On the other hand if $\sigma_a(\varepsilon)$ is small this
means that any departure of $\nu_a$ from its optimal value will
severely worsen the plant performance. Note that this uncertainty is
an estimate of how wide the minimum is. In our one dimensional example
of Fig.~\ref{fig:multi} the uncertainty of the variable $x$ is about
0.7 for the function A and 0.08 for the function B using a value of 
$\varepsilon=0.001$. 

When missing the minimum for a variable $i$ can be compensated by
tuning the variable $j$ we say that the variables $i$ and $j$ are
correlated. They share part of the information of a plant design. For
example the tower height and the aperture inclination must be
correlated, since we need to adjust the aperture inclination to
correct an error in the tower height. We can say that both variables
``share'' the information that the aperture need to look directly to
the field in order to increase the interception efficiency. This
correlation can be 
mathematically expressed by the correlation coefficient between the
two variables $\rho_{ij}$, a real
number between -1 and 1 that measures the difference between the
quantity $\sigma_i(\varepsilon)$ and its would-be value if the
variable $j$ would have stayed fixed at its minimum value $\bar
\nu_J$ (the quantity $\sigma_{\rm inner}$ in
Fig.~\ref{fig:errorelipse}). 

A value of $\rho_{ij}$ close to zero means that the variables $\nu_i$
and $\nu_j$ encode different information of the plant design. On the
other hand if $\rho_{ij}$ is close to its extreme values
$\rho_{ij}\approx \pm 1$, both
variables $\nu_i$ and $\nu_j$ encode the same information, since a non
optimal value of one of them can always be compensated by choosing
the other appropriately. 

On one hand if a plant design variable has a very large uncertainty this
variable can be dropped from the plant design, since it plays no role
in the plant performance. On the other hand, ideal plant design
variables should be \emph{as 
  uncorrelated as possible}, since otherwise the optimization
algorithm will be looking for an optimal value of a function with an
almost flat direction. 

We are interested in computing the quantities $\sigma_i(\varepsilon)$
and the correlation matrix $\rho_{ij}$. As it is shown in
the~\ref{ap:hess}, when $\varepsilon$ is small, these quantities can
be estimated by computing the Hessian of the objective function
evaluated at the minimum 
\begin{equation}
  H_{ij} = \left. \frac{\partial^2 E(\nu_i;p_\alpha,I)} 
    {\partial \nu_i\partial \nu_j}\right|_{\nu = \bar \nu}.
\end{equation}
 Being concise, the
quantities $\sigma_i(\varepsilon)$ are given by \emph{the diagonal
  elements of the inverse of the Hessian}
\begin{equation}
  \sigma_i^2(\varepsilon) = (H^{-1})_{ii}\, \varepsilon,
\end{equation}
and the correlation coefficients are given by the normalized non
diagonal elements of the inverse hessian
\begin{equation}
  \rho_{ij} = \frac{(H^{-1})_{ij}} {\sqrt{(H^{-1})_{ii}(H^{-1})_{jj}}} .
\end{equation}
We emphasize here that this expressions are accurate only in the case
that $\varepsilon$ is small enough so that all the
$\sigma_i(\varepsilon)$ are small\footnote{In the almost trivial case that
the objective function is quadratic the above mentioned expressions
are exact for any value of $\varepsilon$, but this is hardly an
interesting case of study.}. For practical purposes we can only
consider the resulting $\sigma_i(\varepsilon)$ an estimate of the real
uncertainties, but this information is more than enough to obtain a
qualitative understanding of how the different design variables impact
the plant performance and to interpret the result of the optimization
that different algorithms produce. 

The strategy that we follow once the optimal value of the design
variables is found is to compute numerically the Hessian matrix, and
then invert it. To compute the numerical first and second derivatives
we use standard finite differences expressions (we drop the non
variable arguments of $E$ for clarity):
\begin{eqnarray}
  \frac{\partial E(\nu_i)}{\partial \nu_k} &=& \lim_{h\rightarrow
    0}\frac{E(\nu_i + h_k\hat k) -  E(\nu_i - h_k\hat k)}{2h}\\ 
  \frac{\partial^2 E(\nu_i)}{\partial^2 \nu_k} &=& \lim_{h\rightarrow
    0}\frac{E(\nu_i + h_k\hat k) - 2E(\nu_i) + E(\nu_i - h_k\hat k)}{h^2}\\ 
  \nonumber
  \frac{\partial^2 E(\nu_i)}{\partial \nu_k\partial \nu_j} &=&
  \lim_{h_l\rightarrow 
    0}\frac{1}{4h_k h_j} \left\{E(\nu_i + h_k\hat k + h_j\hat
    j)\right. \\
  \nonumber
  &-& E(\nu_i
  + h_k\hat k - h_j\hat j) \\
  \nonumber
  &-& 
  E(\nu_i - h_k\hat k + h_j\hat j) \\
  &+& \left. E(\nu_i - h_k\hat k + h_j\hat
    j)\right\}
\end{eqnarray}
where $\hat i$ represents the unit vector in the direction of the
variable $\nu_i$.

Choosing an appropriate step size (the values for $h_i$ in the
previous formulas) is crucial to compute derivatives numerically. This is a
well known problem: numerical differentiation is an ill defined
problem for finite precision arithmetic.  

We have chosen a step size $h_i$ for each direction that meets
the following  two criteria:
\begin{enumerate}
\item The numerical computation of the first derivative gives a very
  small result.
\item The numerical computation of the second derivative gives
  approximately the same result if we vary the step by a 10-20\%. 
\end{enumerate}
The first condition ensures that we are not having rounding errors
(being at the minimum means that the first derivative should be
zero). The second condition ensures that the computation of the second
derivatives is stable. 

\section{Results}
\label{sec:res}

Following the steps described in the previous sections we will present
some results about the optimal field design of some typical solar
power plants. 

For the case of a cavity receiver we will analyze the case of a 900
heliostat power plant. This plant produces aound $15$ MWe. We will
label this plant as \texttt{N900}.  

For the case of the cylindrical receiver we will analyze the case of a
3000 heliostat field. This circular filed layout produces around $50$
MWe. We will label this plant \texttt{C3000}.  

In both cases the optimization criteria consist in finding
the cheapest price for the generated power. 

\subsection{Algorithm analysis}

First we will focus in the \texttt{N900} plant. In Tab.~\ref{tab:min}
we can see the result of the optimization 
procedure with each of the algorithms.  As can be seen the three
algorithm give a similar result for the value of the objective
function at the minimum (quantity $\overline E$). Nevertheless, this
optimal value is achieved with different values of some of the
variables. The value of $a_0$ is 5.4 for the minimum found by the
\texttt{NSPOC} algorithm and $2.3$ for the Genetic one. Since the
value of the objective function is almost the same. This may induce to
think that the objective function has several almost degenerate 
local minimums. 

\begin{table}
  \centering
  \begin{tabular}{|l|l|l|l|l|}
    \hline
    & \texttt{NSPOC} & Genetic & \texttt{MINUIT} &
    $\sigma(\varepsilon=0.001)$ \\
    \hline
    \hline
    $\overline E$  & 0.18092 & 0.18113 & 0.18311 & N/A\\
    \hline
    \hline
    $a_0$      & 5.4   & 2.3   & 3.5 & 2.8\\
    \hline
    $a_1\times 10^2$      &  3.15  &  3.78  & 3.37 & 0.55 \\
    \hline
    $d_\theta$  & -9.5  & -1.9  & -5.4 &  8.1\\
    \hline
    $e_\theta$    &  -0.05  &  -0.78  & -0.01 & 0.66 \\
    \hline
    $\epsilon$      &  0.88  &  0.24  & 1.8 & 1.1\\
    \hline
    $\delta$ &  0.169  &  0.024  & 0.015 & 0.097\\
    \hline
    $b\times 10^3$      & 54.2  & 32.8  & 29.6 & 24.0\\
    \hline
    $D_0^0$      & 16.8  & 19.5  & 19.1 & 2.5\\
    \hline
    $h_T$      &  120.     & 117.    & 123.  & 12.\\
    \hline        
    $r$        & 10.78     & 10.74   & 10.85  & 0.99\\
    \hline
    $e_L$      &  28.6     & 26.5    & 38.1  & 9.3\\
    \hline
  \end{tabular}
  \caption{Comparison of the three algorithms for the \texttt{N900} plant.
    The second, third and fourth
    column shows the results for the \texttt{NSPOC}, the Genetic and the
    \texttt{MINUIT} algorithm respectively. Finally the last column quote
    the values of $\sigma(\varepsilon)$ for a change in the function
    of $\varepsilon=0.001$. The first row show the value of the
    function at the minima. Subsequent rows show the value of each of
    the design variables. . The following lines shows the value of different variables
    at the minima for each algorithm.}
  \label{tab:min}
\end{table}

But nothing could be more wrong. As our minimum analysis shows all
these values corresponds (approximately) to one and the same
minimum. But this minimum is very wide in some of the directions. This
can be seen by computing the uncertainties needed to change the
objective function value in an amount 
$\varepsilon$ for each of the directions. We will choose
$\varepsilon=0.001$ as an estimate of the difference of the optimum
value obtained by different algorithms. Clearly this small change of
the objective function does not change the plant properties,
since in the plant evaluation one uses approximations that lead to an
estimate of the plant performance that is accurate with less precision
than this value of $0.001$. 

As can be seen in Tab.\ref{tab:min} the variables that show a strong 
discrepancy between the results of different algorithms corresponds to
directions in the objective function in which the minimum is very
wide. For the case of $a_0$ the associated uncertainty is 2.8, meaning
that the values 5.42 and 2.27 that different algorithms find
correspond to the same optimum. This confirms that all algorithms are
finding the same minimum, 
but that in this minimum the objective function has a very mild
dependence with some variables. Variables that are crucial for the
plant design are the ones that have 
a small associated uncertainty, like for example $a_1$ or $r$. For
this variables we can see that all algorithms find approximately the
same value of the objective function.

In one sentence, the discrepancy between the value of the value of the
variables at the minimum obtained with different algorithms is nicely
explained once ones analyzes how strong is the dependence of the
objective function at the minimum with respect to different
variables. One can conclude that all the algorithms are finding the
same basic plant design.

Also we have checked that the correlations between parameters are not
large. The largest correlation turns out to be between variables $b_0$ and
$b_1$, and amounts to 0.91. Note that precisely these are the variables
with a large uncertainty. Variables with small associated uncertainties
have usually small correlations. For example in the case of $r$ and
$a_1$ the correlation amounts to $-0.08$. 

The correlation matrix also seem to pick up important information of
the plant design. For example the tower height $h_T$ is positively
correlated with the aperture inclination $e_T$ (the correlation
amounts to $0.1$), indicating that if ones makes the tower higher one
needs to incline more the aperture. 

\begin{table}
  \centering
  \begin{tabular}{|l|l|l|l|l|}
    \hline
    & \texttt{NSPOC} & Genetic & \texttt{MINUIT} &
    $\sigma(\varepsilon=0.001)$ \\
    \hline
    \hline
    $\overline E$  & 0.17022 & 0.16981 & 0.17023 & N/A\\
    \hline
    \hline
    $a_0$      & 0.97  & 0.27 & 0.25 & 0.26\\
    \hline
    $a_1\times 10^2$      & 3.40  & 3.52 & 3.29 & 0.35 \\
    \hline
    $d_\theta$ & -1.4 & -2.8 & -5.3 &  1.9\\
    \hline
    $e_\theta$   & -0.5 & 0.0 & 1.8 & 1.5 \\
    \hline
    $\epsilon$      & 0.29 & 0.37 & 0.16 & 0.55\\
    \hline
    $\delta$ & 0.109 & 0.081 & 0.134 & 0.039\\
    \hline
    $b\times 10^3$      & 40. & 60. & 47. & 13. \\
    \hline
    $D_0^0$      & 17.8 & 15.3 & 15.1 & 1.9 \\
    \hline
    $h_T$      &  145.     & 147.    & 152.  & 13.\\
    \hline        
    $r$        & 8.58   & 8.66  & 8.57 & 0.96\\
    \hline
    $h_r$      &  8.14     & 8.15 & 8.17  & 0.96\\
    \hline
  \end{tabular}
  \caption{Comparison of the three algorithms for the \texttt{C3000} plant.
    The second, third and fourth
    column shows the results for the \texttt{NSPOC}, the Genetic and the
    \texttt{MINUIT} algorithm respectively. Finally the last column quote
    the values of $\sigma(\varepsilon)$ for a change in the function
    of $\varepsilon=0.001$. The first row show the value of the
    function at the minimum ($\overline E$). Subsequent rows show the
    value of each of 
    the design variables. The following lines shows the value of
    different variables at the minimum for each algorithm.}
  \label{tab:min3000}
\end{table}

In Tab.\ref{tab:min3000} we can see the same information for the case
of the \texttt{C3000} plant. As we can see the conclusions are
roughly the same. The three algorithms seem to find the same design as
can be seen by comparing the values of the variables at the
minimum found by different algorithms with their respective
uncertainties. 

Comparing Tab.~\ref{tab:min} and Tab.~\ref{tab:min3000} we observe
that in general the uncertainties are smaller for the \texttt{C3000}
plant. This means that for the design of bigger plants variables
need an accurate value if we want to obtain a high performance. The
design of 
solar power plants become more involved and difficult when the size of
the power plant increases. This conclusion seem intuitively correct,
since there are some variables (ej. $d_\theta$ that ``compress'' the field
layout in the south) that seem to be crucial only for big circular
plant designs, and more or less irrelevant for the design of small
power plants.  

\subsection{Performance analysis}

A summary of the performance of the different algorithms is presented
in Tab.~\ref{tab:time}. 
\begin{table}[h!]
  \centering
  \begin{tabular}{|l|l|l|l|l|}
    \hline
    PLANT & \textbf{Algorithm}& Calls & Time $[s]$ 
    & Time $[a.u.]$\\
    \hline
    \multirow{3}{*}{\texttt{N900}} 
    &\texttt{NSPOC}  & 1270 & 17280 & 1.00 \\
    &\texttt{MINUIT} & 3737 & 86450 & 5.00 \\
    &Genetic         & 2381 & 40539 & 2.35 \\
    \hline
    \multirow{3}{*}{\texttt{C3000}} 
    &\texttt{NSPOC}  & 1872  &  41460 & 1.00 \\
    &\texttt{MINUIT} & 4599  & 102699 & 2.47 \\
    &Genetic         & 10398 & 254829 & 6.14 \\
    \hline
  \end{tabular}
  \caption{Running times of the different algorithms in our prototype
    plants. The first three columns corresponds to the 900 heliostat
    north field plant \texttt{N900}, and the last three columns to
    the 3000 heliostat cylindrical receiver plant \texttt{C3000}. The
    table shows the number of calls to the plant evaluation function as
    well as the time of the run, both in seconds and relative to the
    running time of the \texttt{NSPOC} algorithm. As the reader can
    see, it took for the global algorithms (both \texttt{MINUIT} and
    Genetic) between 2.5 and 5 times more time to find the optimum
    plant design than the time used by \texttt{NSPOC} for the case of
    the plant \texttt{N900}. In the case of the plant \texttt{C3000}
    the case is worse (between 4 and 7 times more time).}
  \label{tab:time}
\end{table}

The \texttt{NSPOC} algorithm is always faster than the global
optimizers \texttt{MINUIT} and Genetic. In the case of the plant
\texttt{N900}, \texttt{MINUIT} uses 5 times more time to achieve the
optimum, while Genetic uses 2.35 more time than \texttt{NSPOC}. This
difference is similar for the case of the \texttt{C3000}
plant. 

The two global optimizers that we have used are stochastic in nature,
thus this running times should not be treated as exact numbers, but
they are representative, and the conclusion is always the same: the
\texttt{NSPOC} algorithm outperforms the global optimizers while
giving the same results.  

As we have said, the \texttt{NSPOC} code only re-compute the
``needed'' pieces of the objective function. For example blocking and
shadowing effects are not computed when we perform a line search in the
direction of the aperture size. This approach seems
to be very successful for the case of the \texttt{N900} plant. On
average each function call takes between a 25\% and a 
50\% less time for the \texttt{NSPOC} algorithm than for the
others. In the more involved case of the \texttt{C3000} plant the
\texttt{NSPOC} function calls still are faster, but with a large
margin. 

\section{Conclusions and perspectives}
\label{sec:conc}

We have presented and analyzed a method to design solar power
plants. Based in the field layout done within the
German-Spanish GAST project (1982-86) in the company INTERATOM, we
have proposed a method to design solar power plants. This method
reduces the plant design to the value of 11 variables that determine
the field layout, tower and receiver characteristics. The problem of
finding optimal plant designs is reduced to the numerical problem of
optimizing a non linear function.

Due to the non linear nature of the target function and the high
computational costs of evaluating plant performances, it is crucial 
to find a robust and fast algorithm to perform plant designs. Fast
algorithms to optimize functions of several variables are local in
nature: they find a local optima of the target function, but give no
information on the possible existence of other global optima. On the
other hand global optimizers are much more computationally expensive
and their coding is far more involved. 

In this work we have tackled the optimization problem with different
algorithms, both local and global. We have shown that in all the cases
our local NSPOC algorithm give the same results as other more complex
optimizers with significantly less
computing effort. 

The design variables that we have choosen have enough flexibility
to provide solar plant layouts with very high performance for a broad
range of sizes: from small $10$ MWe plants to big $200$ MWe
plants. With our choice of variables we find that our local optimizer
outperforms global optimizers.

In analyzing the result of different algorithms we have also developed
an interesting method to get information on how different variables
affect the performance of an optimal design. The method is based in
analyzing the Hessian of the objective function at its optimum value,
and gives us an estimate of how much the departure of a variable from
its optimal value affects the plant performance. We have observed that
circular-like big plants need a more accurate tuning of the variables
in order to achieve an optimal performance. 

This information can be used to speed up the optimization
process. It turns out to be convenient to first tune variables that are
crucial for the optimal design, and only after worry about variables
that have a mild impact in the plant performance. Our \texttt{NSPOC}
algorithm partially profit from this information to speed up the
convergence to the optimum. 

Moreover, this method of evaluating the impact of design variables in
the plant performance is key to the study of the improvement of plant
designs. Any new proposed design \emph{should} include a similar
analysis to detect superfluous variables in the design and determine
what are the key ingredients of the new design. 

This works provide the tools to address very interesting
problems, what we call \emph{parametric analysis}, in which the
efficiency of the solar power plant can be studied as a function of
plant design parameters. Questions like Does a terrain with slope
improve the plant efficiency? How much does the heliostat aspect ratio
affect the plant efficiency? Can we build cheaper plants by using
cheaper heliostats without losing efficiency? Some of these questions
are currently under study by the authors, and preliminary results were
presented in the solarPACES 2011 conference~\citep{Ramos:2011aa}.

We have also provided the basis to analyze more complex plant
designs. An example that we consider very interesting are multi-tower
layouts, in which a solar power plant is built with several towers,
and heliostats choose which tower they aim based on performance (some
results were already presented in the 2009 solarPACES
conference~\citep{Ramos:2009aa}), but there are many more
possibilities like multi-cavity receivers, or properly addressing the
scalability problem in solar power plant designs: How can we build
small plants that can later be enlarged without loosing performance?
Some of this questions and others within the framework presented here
are currently being studied. 

\appendix

\section{Computation of uncertainties and correlations.}
\label{ap:hess}

The techniques described in this appendix are typical in statistical
description of data. In this context the function that one wants to
optimize is the quadratic deviation between data and predictions of a
model (usually called $\chi^2$). Thus any interested reader can
consult any standard book on statistics for a more detailed proof of
the expressions developed here. We will assume that we are
\emph{minimizing} an objective function. The expressions remain
basically unchanged for the case of the maximization of a function if
one changes the sign of $M$, $(M^{-1})$ and the Hessian in the
expressions that follow.  

We will start assuming that our objective function of $n$ variables
($x_i, i=1\dots,n$)  is 
quadratic with a minimum located at $x_i=\bar x_i$. The most general
function with these characteristics can be written as
\begin{equation}
  f(x) = f_0 + \sum_{ij}(x_i - \bar x_i) M_{ij} (x_j -
  \bar x_j)
\end{equation}
where $f_0$ is the value of the function at the minimum and $M_{ij}$
is a symmetric positive definite matrix. To obtain  
$\sigma_a(\varepsilon)$ we will define a new function 
$g(x_k)$ ($k\neq a$) of $n-1$ variables (all but the $a-th$) equal to
the value of the original function with $x_a$ 
fixed at $x_a = \bar x_a + \sigma_a(\varepsilon)$.
\begin{eqnarray}
  \nonumber
  g(x_k) &=& f\left(x_k,x_a=\bar
    x_a +\sigma_a(\varepsilon)\right) \\
  \nonumber
  &=& f_0 + \sum_{i,j\neq a}(x_i - \bar x_i) M_{ij} (x_j -
  \bar x_j) \\
  &+& 2 \sum_{i\neq a}(x_i - \bar x_i) M_{ia} \sigma_a(\varepsilon)
  + M_{aa}\sigma_a^2(\varepsilon).
\end{eqnarray}
Taking the gradient of $g(x_k)$ and equating it to zero we obtain
the position of the minimum value of the objective function at fixed
$x_a = \bar x_a + \sigma_a(\varepsilon)$. 
\begin{equation}
  \frac{\partial g}{\partial x_k} = 0 \Longrightarrow
  x_k^* = \sum_{i\neq a}\left(m^{-1}\right)_{ki} M_{ia} \sigma_a(\varepsilon)
\end{equation}
where $m$ is the minor matrix that corresponds to the element $(a,a)$
(this is the matrix that results from $M$ cutting down the $a^{th}$
column and rows), and $m^{-1}$ is its inverse. The value of the
objective function at this point is given by 
\begin{eqnarray}
  \nonumber
  g(x_k^*) &=& f\left(x_k^*,x_a=\bar
    x_a +\sigma_a(\varepsilon)\right) \\ 
  &=& f_0 + \left[
    M_{aa} - \sum_{i,j\neq a} M_{ai}\left( m^{-1}\right)_{ij}M_{ja} 
  \right]\sigma_a^2(\varepsilon)
\end{eqnarray}
the quantity between brackets can easily be recognised as
$1/\left(M^{-1}\right)_{aa}$ (the inverse of the diagonal entry of the
original matrix $M$). The condition that determines
$\sigma_a(\varepsilon)$ is that the increase in the function respect
the value at the minimum should be $\varepsilon$. So the quantity
$\sigma_a(\varepsilon)$ is given by
\begin{equation}
  \sigma_a^2(\varepsilon) = {\left(M^{-1}\right)_{aa}}\, {\varepsilon}.
\end{equation}

If we define $\sigma_{{\rm inner}}(\varepsilon)$ as the would-be 
$\sigma_a(\varepsilon)$ without tuning the variable $j$, it is clear
that a comparison between $\sigma_a(\varepsilon)$ and 
$\sigma_{{\rm inner}}(\varepsilon)$ would give us information about
the correlation between variables $i$ and $j$. If we define the 
correlation coefficient between variables $i$ and $j$ by the equation 
\begin{equation}
  \rho_{ij} = \frac{(M^{-1})_{ij}} {\sqrt{(M^{-1})_{ii}(M^{-1})_{jj}}}
\end{equation}
we can easily check that $\sigma_{{\rm inner}}(\varepsilon)$ and
$\sigma_a(\varepsilon)$ are related by
\begin{equation}
  \sigma^2_{\rm inner} = (1-\rho^2_{ij})\sigma^2_j(\varepsilon).
\end{equation}
Note that in the case $\rho_{ij} = 0$ the two quantities are equal. 

In the case that the objective function is not quadratic, we can use
Taylor theorem. Close enough to the minimum \emph{any} function is
well approximated by a quadratic function. This means that for a
general function the role of $M_{ij}$ is played by the Hessian
evaluated at the minimum 
\begin{equation}
  M_{ij} = \left.\frac{\partial^2
      f(x_i)}
    {\partial x_i\partial x_j}\right|_{x=\bar x} + \mathcal
  O(x-\bar x)^2.
\end{equation}
It is worth noting that this approximation will fail if
$\sigma_a(\varepsilon)$ is large. In order to ensure that this does
not happen one need to keep $\varepsilon$ sufficiently small.


\end{document}